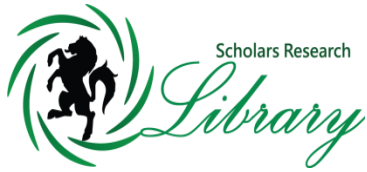

# Scholars Research Library

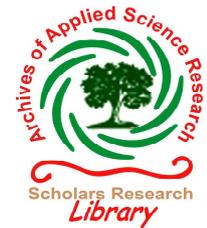

## A Minimum Angle Method for Dual Feasibility

**Syed Inayatullah, Nasiruddin Khan, Muhammad Imtiaz, Fozia Hanif Khan**

*Department of Mathematics, University of Karachi, Karachi, Pakistan, 75270*

______________________________________________________________________

**Abstract**

In this paper we presented a new driving variable approach in minimum angle rule which is simple and comparatively fast for providing a dual feasible basis. We also present experimental results that compare the speed of the minimum angle rule to the classical methods. The experimental results showed that this algorithm outperforms all the previous pivot rules.

______________________________________________________________________

## Introduction

Linear programming has been a fundamental topic in the development of the computational sciences [6]. Algorithms for finding a solution to a system of linear inequalities first studied by Fourier in the 19th century [16] and later rediscovered by several mathematicians such as [7] and [5].The revival of interest in the subject in the 1940's was spearheaded by [4] (Simplex Method) in USA.

In general there are two very important and crucial steps involved in almost every variant of simplex method, the selection of entering and leaving basic variable. There are several methods has been developed for this purpose.

Selection of entering basic variable actually is the selection of edge emanating from current extreme point and leaving basic variable is the selection of contrary constraint to that edge. For selection of the entering basic variable there are several rules including Dantzigs largest coefficient rule [9], Largest increase rule, Steepest edge rule [10], Devex rule [11], Largest distance rule [12].

_______________________________________________________________________________

For the selection of leaving basic variable the Dantzig suggests a constraint with minimum right hand ratio. This rule retains the primal feasibility through out the process but highly affected by the degeneracy.

In the current paper we shell use the idea behind the largest coefficient rule and the steepest edge selection rule for entering basic variable in the selection of leaving basic variable. Both of these rules select the edge of movement, having maximum projection of objective gradient. Now we extend this idea for the selection of leaving basic variable. Actual reasonable and appealing logic behind the largest coefficient and especially steepest edge rule is that the edge more likely to reach the near optimal extreme point is the edge having largest projection of objective vector.  A similar idea for leaving basic variable have been suggested by [13],

"The constraint with a larger cosine value with objective gradient is more likely to be binding at an optimal extreme point than one with a smaller value."

[17] presents ratio test free rule based on a similar idea describing a similar description about optimal solution. Pan claimed that the most obtuse angle ratio test free rules are best possible rules for the selection of leaving basic variables.

A related and seemingly good result by [14] and [15] "The non-redundant constraint which is more parallel to the gradient of objective function is binding at optimal"

On the basis of above assertions about constraints binding at optimal point, here in this paper we exhibit a simple and fast algorithm for solving linear programs. The algorithm consists of two phases First phase obtains the dual feasibility by using minimum angle rule and in the second phase we use the dual simplex method to obtain the optimality. The advantage of new approach is that it could start with a primal infeasible basic solution, so introduction of artificial variables is not mandatory. Additionally, in this paper we have compared its efficiency to the classical finite pivot methods for attaining dual feasibility in random linear programs.

**The Basic Notations**

A general LP problem is a collection of one or more linear inequations and a linear objective function involving the same set of variables, say $x_1,…,x_n$, such that

$$\text{Maximize} \quad c^T x$$

(1)     subject to     $Ax = b$ ,
                        $x \geq 0, x \in \Re^n$

where $x$ is the decision variable vector, $A \in \Re^{m \times n}$, $b \in \Re^m$, $c \in \Re^n$.

We shell use $a_i$ and $a_{.j}$ to denote the $i^{th}$ row and $j^{th}$ column vector of $A$. Now we define a basis $B$ as an index set with the property that $B \subseteq \{1,2,...n\}$, $/B/= m$ and $A_B := [a_{.j} / j \in B]$ is an invertible matrix, and Non-basis $N$ is the complement of $B$. i.e. $N := \{1,2,...n\}\backslash B$.

Given a basis $B$ and $x=(x_B, x_N)$, we can have an equivalent canonical form of above system (1),

$$\text{Maximize} \quad \overline{c}^T x_N + z$$

(2) $\quad$ subject to $\quad \overline{A} x_N + x_B = \overline{b},$
$$x_N \geq 0,\ x_B \geq 0,$$
$$x_N \in \Re^N,\ x_B \in \Re^B.$$

So, we may construct the following dictionary for basis $B$ c.f. [3],

$$D(B) = \begin{bmatrix} z & -\overline{c}^T \\ \overline{b} & \overline{A} \end{bmatrix}$$

where,
$$\overline{A} = A_B^{-1} A_N,$$
$$\overline{c}^T = c_N^T - c_B^T A_B^{-1} A_N,$$
$$\overline{b} = A_B^{-1} b$$
$$z = C_B^T A_B^{-1} b$$

The associated basic solution could directly be obtained by setting $x_B = \overline{b}$. Here onward in this text we assume that the reader is already familiar about pivot operations, duality and primal-dual relationships. It is well known that if $d_{B0} \geq 0$ then $B$ (or $D(B)$) is called primal feasible (dual optimal); if $d_{0N} \geq 0$ then $B$ (or $D(B)$) is called primal optimal (dual feasible); and if both $d_{B0} \geq 0$ and $d_{0N} \geq 0$ then $B$ (or $D(B)$) is called optimal feasible. A basis $B$ (or a dictionary) is called inconsistent if there exists $i \in B$ such that $d_{i0} < 0$ and $d_{iN} \geq 0$, and unbounded if there exists $j \in N$ such that $d_{0j} < 0$ and $d_{Bj} \leq 0$.

**Half space representation of an LP:**

For every canonical form of an LP we can have an inequality form (also called half space representation) by removing basic variables from the constraints. The resulting inequality system represents optimization of objective vector over a polyhedron (intersection of half spaces) in non-basic space. for example system (2) has the following Half space representation,

$$\text{Maximize} \quad \overline{c}^T x_N + z$$

(3) $\quad$ subject to $\quad \overline{A} x_N \leq \overline{b},$
$$x_N \geq 0,\ x_N \in \Re^N$$

In this inequality form the vector coefficient of each constraint will be called $h$-vector and the directions along non-basic variables $x_N$ will be called $e$-vectors. A notation of $h(B,i)$, $i \in B$ will be used for $h$-vector of the constraint corresponding to basic variable $x_i$, and $e(B,j)$, $j \in N$. will be used for $e$-vector along the non-basic variable $x_j$. The $e$-vectors along which the objective value is increasing are called improving edges. For any improving edge we call a constraint as contrary if it restricts the improvements in the objective function along that direction.

**Lemma**

If $i \in B$ and $J \subseteq N$ then the vector $h(B,i)$ is contrary to $c(B,J)$ if and only if $h(B,i)^{\mathrm{T}} \cdot c(B,J) > 0$.

*Proof*:

Since all the primal constraints are of "$\leq$" form then along the ascent $c(B,J)$ the chance of violation of constraint $C_i$ is valid if and only if $h(B,i)^{\mathrm{T}} \cdot c(B,J) > 0$. ∎

*Remark:* For an LP to be solvable, there must be at least one resisting $h$-vector with respect to each of the improving $e$-vectors.

Consider an improving edge $E_j$ with some resisting constraints and let $K_j$ is index set of all these resisting constraints i.e. $K_j = \{k \mid \overline{a_{kj}} > 0\}$, we call a resisting constraint $C_m$ as *the most contrary*, if $cos\,[\theta(m,j)] \geq cos\,[\theta(k,j)]$, for all $k \in K_j$. Where notation $\theta(a,b)$ denotes angle between $C_a$ and $E_b$.

**Theorem:** Given a solvable LP, If in any associated LP polyhedron there is only one improving $e$dge with only one resisting constraint then that constraint would be the optimal constraint.

*Proof*:-

If an improving edge $E_k$ has only one resisting constraint $C_l$, it means that $C_l$ is the only restriction applied to the objective value along the direction $e(B,k)$. So, removal of $C_l$ will change the feasible region as well as the optimal region (in fact the optimal region would become unbounded). Hence $C_l$ is an optimal (Binding) constraint of the given LP. ∎

On the basis of several assertions described in section 1, we may conveniently conclude the following statement,

**Statement 1**

"Given a solvable LP, if there is only one improving edge and more than one contrary constraint to it then most contrary constraint is more likely to be bind at optimal feasible point."

Just like all the other authors worked on minimum angle rules, we also cannot able to give any theoretical proof of the above statement but here on the basis of this statement we are giving an algorithm **Algorithm 1**. Clear successes of our algorithm on a large variety of random LPs, as shown in section 3, give us evidence about correctness of the statement.

Generally, in an LP there would be more than one improving edge. But we may transform that LP to an equivalent LP having only one improving edge. Here is the method which will give us such transformation.

_______________________________________________________________________________

**SIE Transformation: (Single Improving Edge)**

If a dictionary $D(B)$ has more than one dual infeasible variable (candidate for entering basic variable), geometrically it indicates occurrence of more than one improving edge direction. These directions could be reduced to only one direction by introducing artificially a new non-negative basic variable $x_r$ (called the *driving* variable), equal to following linear combination,

$$x_r = -\sum_{j \in L} d_{0j} x_j \text{ , where } L = \{ j : d_{0j} < 0 \text{ , } j \in N \}.$$

and then by making $x_r$ the leaving basic and any of $\{x_j , j \in L\}$ as entering basic variable (preferably select the variable having most negative coefficient in the objective row). Note that here coefficients of the linear combination are values of dual infeasible variables.

**<u>Problem 1:</u>**

Given a dictionary $D(B)$ of an LP with improving edges $\{E_i : i \in L, L \subseteq N, |L| \neq 1\}$. For an element $l \in L$ transform the LP into an equivalent LP with only one improving edge $E_l$.

**or equivalently,**

Given a dictionary $D(B)$ of an LP with a collection of dual infeasible variables $s_L$ where $L \subseteq N$, and $|L| > 1$. For an element $l$ of $L$ transform $D(B)$ into an equivalent dictionary with only one dual infeasible variable $s_l$.

**Procedure (SIE Transformation)**

**Step 1:** Introduce a new basic variable $x_r$ such that $r := |B| + |N| + 1$ and $x_r := 0$. Set $B := B + \{r\}$.

**Step 2:** Insert the row vector $d_{rj} = .\begin{cases} d_{0j} & , j \in L \\ 0 & , j \notin L \end{cases}$ into the $D(B)$.

**Step 3:** Make a pivot on *(r,l)*. Set $B := B - \{r\} + \{l\}$, $N := N + \{r\} - \{l\}$. **Done.**

**<u>Problem 2:</u>**

Given a dictionary $D(B)$, obtain the dual feasibility.

**Algorithm 1: (Minimum angle rule)**

**Step 1:** If $d_{0j} \geq 0$ for all $j \in N$ then $D(B)$ is dual feasible. **Done.**

**Step 2:** Otherwise let $L$ be a maximal subset of $N$ such that $d_{0j} < 0$, for all $j \in L$.
          i.e. $L = \{ j : d_{0j} < 0 , j \in N \}$

**Step 3:** Let $l$ be an element of $L$ s.t. $d_{0l} \leq d_{0j} \, \forall \, j \in L$. If $|L| = 1$ then set $r := l$ and go to **Step 5.**
          (Ties could be broken arbitrarily)

_______________________________________________________________________________

**Step 4:** Apply SIE Transformation to $D(B)$ taking $e(B,l)$ as the main direction.

**4a:** Introduce a new basic variable $x_r$ such that $r := |B| +|N| +1$ and $x_r := 0$. Set $B := B +\{r\}$.

**4b:** Insert the row vector $d_{rj} = .$
$$\begin{cases} d_{0j} & , j \in L \\ 0 & , j \notin L \end{cases} \text{ into } D(B).$$

**4c:** Make a pivot on $(r,l)$. Set $B := B - \{r\} + \{l\}$, $N := N + \{r\} - \{l\}$.

**Step 5:** Let $K$ be subset of $B$ s.t. $K = \{i : d_{ir} > 0, i \in B\}$. If $K = \varnothing$ then stop. $D(B)$ is dual inconsistent.

**Step 6:** Otherwise, choose $m \in K$, such that $cos[\theta(m,r)] \geq cos[\theta(k,r)] \ \forall \ k \in K$. (Ties could be broken arbitrarily.)

**Step 7:** Make a pivot on $(m,r)$. Set $B := B - \{m\} + \{r\}$, $N := N + \{m\} - \{r\}$.

**Step 8:** Remove $[d_{rj} : j \in N]$ from $D(B)$ and set $B := B - \{r\}$. Go to **Step 1**.

**Explanation:**

Above algorithm is totally developed on the idea revealed from aforementioned statement 1. It ends up on super-optimal (probably optimal) basic solution by successively selecting either optimal or super-optimal constraint at the end of each iteration.

In the algorithm first we check for negative elements in the objective row. Now there are three cases may arise,

1) There is no negative element in objective row,
2) There is only one negative element in the objective row,
3) There are more than one negative element in the objective row.

In the first case then end the algorithm with massage "**Done**". In the second case skip **Step 4** and proceed to **Step 5**. In the third case perform **Step 4** (The SIE Transformation) and return a new dictionary having only one negative element (coefficient of a new non-basic variable $x_r$) in the objective row. At this moment we will find set $K$ as mentioned in **Step 5** if $K$ is empty then end the algorithm with massage "**LP is unbounded**" otherwise we will compute $cos [\theta(k,r)] \ \forall \ k \in K$ and determine $m$ such that $cos[\theta(m,r)] \geq cos[\theta(k,r)] \ \forall \ k \in K$. Perform a pivot operation $(m,r)$. Then we may remove the contents of row vector $d_{rN}$ from the dictionary to reduce the size of dictionary, and go to **Step 1** to perform the next iteration.
See also example given at the end of paper.

<u>**Problem 3:**</u> **(Optimality Problem)**

Given an *n*-dimensional vector *c*, *n* x *m* matrix *A*, and an *m*-dimensional vector *b*, solve the linear program max{ $c^T x$ / $x \geq 0$, $Ax \leq b$ }.

<u>**Algorithm 2:**</u>

**Step 1:** Construct the Dictionary $D(B)$ for  an arbitrary basis *B*

**Step 2:** Obtain the dual feasibility by using Algorithm 1.

**Step 3:** If step 2 returns a dual feasible solution, use dual simplex method to obtain optimality.

**3. Experimental Results**

In this section we examined the efficiency of our newly introduced algorithm in contrast to other famous finite pivot rules: Dantzig's largest-coefficient simplex method (with lexicographic minimum ratio test), Bland's smallest index pivot rule, and the ' *b* - rule [1] in the context of attaining dual feasibility. Although these algorithms have a different number of computational Steps per iteration, we focused solely on the number of iterations required to reach dual feasibility. We were not concerned with the amount of computational time required.

Here we present the random models suggested by [2] for generating feasible and infeasible linear programs. The results are obtained by using a C++ compiler under a precision of 7 decimal places.

The suggested random LP model by [2] to test the performance of any dual feasibility attaining algorithm is

$$\text{maximize} \sum_{j=1}^{n} c_j x_j$$

$$\text{subject to} : \sum_{j=1}^{n} a_{ij} x_j \leq b_i \qquad (i = 1,2,...,m)$$

$$x_j \geq 0 \qquad\qquad (j = 1,2,...,n)$$

We generated 500 linear programs with the coefficients $c_j$,  $b_i$, and  $a_{ij}$ chosen randomly from the integer interval [-50, 50].We present the results of these tests in the following table.

_______________________________________________________________________________

### Table 1: Comparison with classical methods on low and high dimensional LP.s

|        | B' Rule   | Bland's rule | Dantzig  | Minimum angle |
|--------|-----------|--------------|----------|---------------|
| 3x3    | 1.37 (1.02) | 2.07 (1.24) | 2.06 (1.23) | 1.44 (0.82) |
| 3x5    | 1.76 (1.04) | 2.42 (1.10) | 2.36 (1.08) | 1.91 (0.84) |
| 3x7    | 1.76 (1.35) | 2.44 (1.36) | 2.33 (1.22) | 2.12 (0.85) |
| 5x5    | 3.14 (1.98) | 3.60 (1.56) | 3.21 (1.32) | 2.26 (1.05) |
| 5x10   | 4.07 (2.32) | 5.09 (3.51) | 3.66 (1.33) | 2.82 (1.22) |
| 7x 5   | 4.43 (2.45) | 4.62 (1.77) | 4.30 (1.52) | 2.63 (1.23) |
| 7x10   | 6.69 (4.03) | 5.79 (2.08) | 4.29 (2.10) | 3.66 (1.52) |
| 10x 5  | 6.69 (4.03) | 5.59 (2.28) | 4.91 (1.88) | 4.23 (1.61) |
| 10x10  | 11.18 (6.13) | 9.13 (3.73) | 6.84 (2.34) | 5.72 (1.84) |
| 10x20  | 14.32 (7.64) | 11.44 (6.00) | 6.66 (2.44) | 5.89 (1.95) |
| 15x15  | 29.26 (13.96) | 19.80 (8.41) | 10.98 (3.82) | 7.98 (2.74) |
| 15x20  | 38.61 (16.96) | 25.57 (10.51) | 12.33 (3.96) | 9.28 (2.86) |
| 20x20  | 55.72 (23.52) | 32.43 (10.89) | 15.93 (4.89) | 12.88 (3.46) |
| 20x30  | 81.34 (36.76) | 43.24 (18.12) | 17.65 (6.80) | 14.11 (3.93) |
| 30x30  | 198.47 (63.07) | 78.25 (20.63) | 30.08 (9.81) | 23.01 (6.57) |
| 40x40  | 467.61 (168.50) | 150.29 (46.28) | 46.32 (13.03) | 32.45 (16.91) |

Our results shows that our new minimum angle pivot rule outperforms almost all the algorithms by far in lower and as well as in higher dimensional linear programs. However, each iteration of our algorithm is more time consuming as compared to others, as in each iteration we may have to compute the dictionary twice. Besides if we compare the results with respect to total number of pivots performed, it would get a very close competition with Dantzig's largest coefficient pivot rule. But in contrast to Dantzig's simplex algorithm, our algorithm does not require any initial basic feasible solution. It means that we may proceed with this method on an LP having neither primal nor dual initial feasible solutions.

**Conclusion**

In this paper we developed the minimum angle algorithm, for the solution of general linear programs. At the end we also showed (through experimental results on random LPs) that our new algorithm beats all the classic ones.

**Example: (Dual Feasibility)**

Consider the following dictionary with $B = \{3, 4, 5\}$ and $N = \{1, 2\}$

|   |    | 1  | 2  |
|---|----|----|----|
|   | 0  | -3 | -5 |
| 3 | 4  | 1  | 0  |
| 4 | 12 | 0  | 2  |
| 5 | 18 | 3  | 2  |

**Iteration 1:** Here $L = \{1,2\}$. We chose $l = 2$, and after SIE Transformation taking $r = 6$, we get $K = \{4,5\}$.

|   |    | 1  | 2  |
|---|----|----|----|
|   | 0  | -3 | -5 |
| 3 | 4  | 1  | 0  |
| 4 | 12 | 0  | 2  |
| 5 | 18 | 3  | 2  |
| 6 | 0  | -3 | -5 |

set $B := \{3,5,2\}$.

|   |    | 1    | 6    | $\cos[\theta(k,6)]$ $k \in K$ |
|---|----|------|------|------|
|   | 0  | 0    | -1   |      |
| 3 | 4  | 1    | 0    |      |
| 4 | 12 | -1.2 | 0.4  | 0.316 ←max |
| 5 | 18 | 1.8  | 0.4  | 0.217 |
| 2 | 0  | 0.6  | -0.2 |      |

|   |    | 1  | 4   |
|---|----|----|-----|
|   | 30 | -3 | 2.5 |
| 3 | 4  | 1  | 0   |
| 6 | 30 | -3 | 2.5 |
| 5 | 6  | 3  | -1  |
| 2 | 6  | 0  | 0.5 |

Pivot on cell (4, 6), delete $[a_{6j}: j \in N]$ and

|   |    | 1  | 4   |
|---|----|----|-----|
|   | 30 | -3 | 2.5 |
| 3 | 4  | 1  | 0   |
| 5 | 6  | 3  | -1  |
| 2 | 6  | 0  | 0.5 |

→     → →

**Iteration 2:** Here $L = \{1\}$, $K = \{3,5\}$. Pivot on (3,1).

➜

**Dual Feasible**

|   |    | **1** | **4** | Cos[θ(k,1)] |
|---|----|----|-----|--------------|
|   | 30 | -3 | 2.5 | $k \in K$ |
| **3** | 4 | 1 | 0 | 1 ← max |
| **5** | 6 | 3 | -1 | 0.949 |
| **2** | 6 | 0 | 0.5 | 0 |

|   |    | **3** | **4** |
|---|----|----|-----|
|   | 42 | 3 | 2.5 |
| **1** | 4 | 1 | 0 |
| **5** | -6 | -3 | -1 |
| **2** | 6 | 0 | 0.5 |